\documentclass[twocolumn]{article}

\pdfoutput=1

\usepackage[T1]{fontenc}
\usepackage[utf8]{inputenc}
\usepackage[english]{babel}
\usepackage{graphics}
\usepackage{epsfig}
\usepackage{amsmath} 
\usepackage{amssymb}
\usepackage{bbold}
\usepackage{url}
\usepackage{hyperref}
\usepackage{float}
\usepackage{pifont}
\usepackage[table]{xcolor}

\usepackage{draftwatermark} 
\SetWatermarkText{Preprint} 
\SetWatermarkScale{5} 

\definecolor{lightgray}{gray}{0.8}

\DeclareMathOperator{\R}{\mathbb{R}}					
\DeclareMathOperator{\1}{\mathbb{1}}					

\title{Model Reduction for Complex Hyperbolic Networks}
\author{Christian Himpe\thanks{Contact: \href{mailto:christian.himpe@uni-muenster.de.de}{\nolinkurl{christian.himpe@uni-muenster.de}}, \href{mailto:mario.ohlberger@uni-muenster.de}{\nolinkurl{mario.ohlberger@uni-muenster.de}}, Institute for Computational and Applied Mathematics at the University of M\"unster, Einsteinstrasse 62, D-48149 M\"unster, Germany} \and Mario Ohlberger\footnotemark[1]{}}
\date{}

\begin{document}

\setlength{\parindent}{0pt}

\maketitle
\thispagestyle{empty}
\pagestyle{empty}


\begin{abstract}\bf
We recently introduced the joint gramian for combined state and parameter reduction [C.~Himpe and M.~Ohlberger. Cross-Gramian Based Combined State and Parameter Reduction for Large-Scale Control Systems. arXiv:1302.0634, 2013], which is applied in this work to reduce a parametrized linear time-varying control system modeling a hyperbolic network.
The reduction encompasses the dimension of nodes and parameters of the underlying control system.
Networks with a hyperbolic structure have many applications as models for large-scale systems.
A prominent example is the brain, for which a network structure of the various regions is often assumed to model propagation of information.
Networks with many nodes, and parametrized, uncertain or even unknown connectivity require many and individually computationally costly simulations.
The presented model order reduction enables vast simulations of surrogate networks exhibiting almost the same dynamics with a small error compared to full order model.
\end{abstract}

~


\section{Introduction}
Complex Networks are often employed as models for large-scale systems like connectivity inside the brain, linking structure of the Internet or trust relations in social networks.
Even in cosmology, causality can be modeled based on a network as demonstrated in \cite{krioukov12}.
Such networks are of hyperbolic structure, in which older nodes are favorably connected compared to younger ones.
In many settings the interconnections of a network contain uncertainties or represent the (possibly unknown) parametrized quantities of interest.
Parametrized models with high-dimensional state and parameter spaces are often infeasible to evaluate many times for different locations of the parameter space.
This is due to two effects.
First, the high-dimensional state space makes each integration of the dynamic system computationally costly.
Second, the high-dimensional parameter space may make many simulations necessary.
In this situation model reduction will accelerate these otherwise costly experiments.
Particularly, the combined reduction of state and parameter space will be illustrated.

This setting for model reduction was inspired by \cite{krioukov12}.
The gramian-based (state) reduction approach originates in (approximate) balanced truncation comprehensively described in \cite{antoulas05}.
An alternative computational method for these gramians, based on proper orthogonal decomposition, was introduced in \cite{lall99} under the name \textbf{empirical gramians}.
For the parameter identification and combined state and parameter reduction, the empirical joint gramian from \cite{himpe13a} is utilized.

In the next section the construction of a hyperbolic network is described.
In section~\ref{sr} the state reduction procedure, then in section~\ref{pr} the parameter identification and combined state and parameter reduction is explained.
For an efficient assembly of the required gramians, the empirical cross gramian is presented in section \ref{eg}.
In section~\ref{uq} the usage of empirical gramians in the context of uncertainty quantification is outlined.
Finally, in section~\ref{nr} a sample network is reduced.

\section{Hyperbolic Network}\label{hn}
Generating a hyperbolic network is a dynamic process with a discrete time space.
The following description is taken from \cite{krioukov12} and \cite[Supplementary Notes II. C]{krioukov12}.
At each time step $t_i$ a new node is born by drawing from a uniform random distribution on the circle $\mathbb{S}^1$ yielding a new node\footnote{For unit curvature.} at $\alpha_i \in [0,2\pi]$ and a radius $r_i = 2\ln\frac{i}{v}$ with network degree $v$.
The new node $x_i$ connects to all existing nodes $x_{0 \hdots i-1}$ that satisfy:
\begin{align*}
 r_j + 2 \ln (\pi - | \pi - | \alpha_i - \alpha_j | |) < 2, \; \forall j < i.
\end{align*}

\begin{figure}
 \includegraphics[scale=0.45]{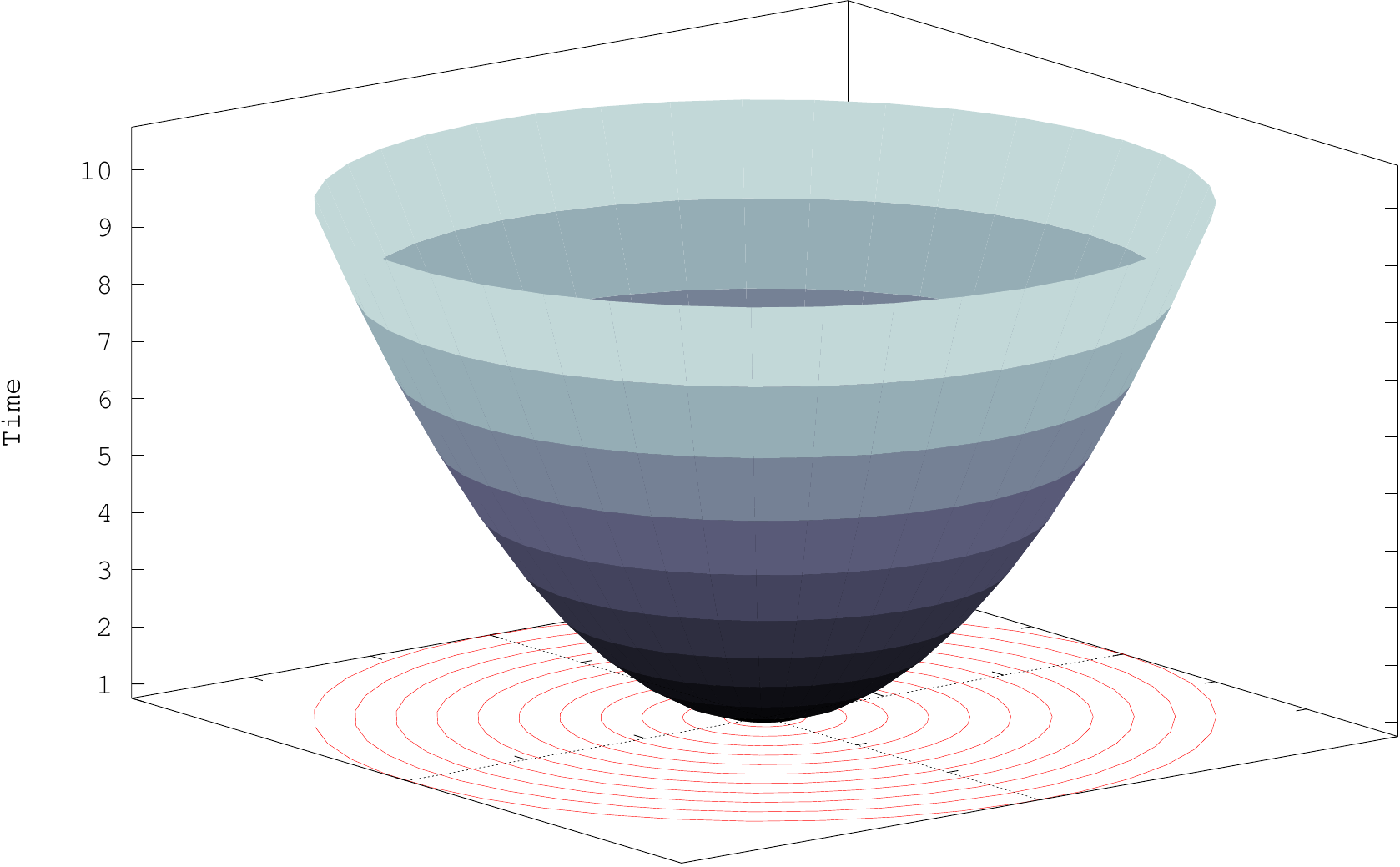}
 \caption{Hyperboloid representing space-time for the hyperbolic network. The contour line represent the projection of the node space at points in time.}
 \label{hy}
\end{figure}
This leads to a "space-time" representation in which the network is created in the shape of a hyperboloid visualized in Figure~\ref{hy}.

In this work a maximum number of $N$ nodes will be set.
Such a hyperbolic network can be modeled as a matrix $A \in \R^{N \times N}$, by treating a matrix element $A_{i,j}$ as connection from node $i$ to $j$, which leads to a dynamic system setting of a basic linear autonomous system:
\begin{align*}
 \dot{x} = Ax,
\end{align*}
with state $x \in \R^N$ and a system matrix $A \in \R^{N \times N}$ embodying the network structure.
Usually some external input or control is applied to the system and the quantities of interest are some subset or linear combination of the systems states.
This leads to a linear control system
\begin{align*}
 \dot{x} &= A x + Bu, \\
       y &= C x,
\end{align*}
with input $u \in \R^{J}$, input matrix $B \in \R^{N \times J}$, outputs $y \in \R^O$ and output matrix $C \in \R^{O \times N}$.
The matrices $B,C$ can for example be used to excite only a certain nodes via $B$ and observe the dynamics in others via $C$.
In this setting it is assumed, the connections, being the components of $A$, between the networks nodes, are parametrized in each component by $\theta \in \R^{N^2}$:
\begin{align*}
 \dot{x} &= A(\theta) x + Bu, \\
       y &= C x.
\end{align*}
As described above, in each time step a new node is born and connects to existing nodes.
Hence, the components of $A$, which are components of $\theta$, too, change over time and yield a parametrized \textbf{linear time-varying control system}:
\begin{align*}
 \dot{x} &= A(\theta(t)) x + Bu, \\
       y &= C x.
\end{align*}
Now, the hyperbolic network can be treated with the model reduction methods of control theory.

\section{State Reduction} \label{sr}
A well known method for model order reduction of the state space in a control system setting is \textbf{balanced truncation}.
In this approach a systems controllability and observability is balanced and the least controllable and observable states are truncated.
A systems controllability is encoded in a gramian matrix $W_C$ which is a solution of the Lyapunov equation $A W_C + W_C A^T = B B^T$.
A systems observability is also encoded in a gramian matrix $W_O$ which is a solution of the Lyapunov equation $A^T W_O + W_O A = C^T C$.
Then, by a balancing transformation, computed from the controllability and observability gramian, the control system is transformed in a manner such that the states are ordered from the most to the least important, for the systems dynamics.
This ordering is based on the Hankel singular values $\sigma_i = \sqrt{\lambda_i(W_O W_C)}$.
After the sorting, the least controllable and observable can be truncated.
The \textbf{cross gramian} encodes controllability and observability into one gramian matrix.
A solution of the Sylvester equation $A W_X + W_X A = -BC$ yields the cross gramian matrix $W_X$ as a solution, if the system is square\footnote{The system has the same number of inputs and outputs}.
In case the system is symmetric, meaning the systems gain $g=-CA^{-1}B$ is symmetric, then the absolute value of the cross gramians eigenvalues equal the Hankel singular values:
\begin{align*}
 |\lambda_i(W_X)| = \sqrt{\lambda_i(W_O W_C)}.
\end{align*}
Given an asymptotically stable system, the cross gramian can also be computed as the time integral over the product of input-to-state and state-to-output map:
\begin{align*}
 W_X = \int_0^\infty e^{At}BC e^{At}dt,
\end{align*}
which will be the basis for computing the empirical gramian variant.

The state reduction is based on the singular values of the cross gramian $W_X$.
A singular value decomposition of $W_X = UDV$, provides a projection of the states $V,U$ in which the states are sorted by their importance.
Without loss of generality, the singular values, composing the diagonal matrix $D$, are assumed to be sorted in descending order.
Based on this projection, the matrices $A,B,C$ and the initial value $x_0$ can be partitioned and reduced,
\begin{align*}
 V &= \begin{pmatrix} V_1 \\ V_2 \end{pmatrix}, \\ 
 U &= \begin{pmatrix} U_1 & U_2 \end{pmatrix}, \\
 \Rightarrow & \begin{cases} \tilde{A_1} = V_1 A U_1 \\ \tilde{B_1} = V_1 B \\ \tilde{C_1} = CU_1 \\ \tilde{x_1}(0) = V_1 x(0) \end{cases}.
\end{align*}
This direct truncation approximates closely the balanced truncation of controllability and observability gramians, but does not require an additional balancing transformation.

In case the system is not square or not symmetric, following the approach from \cite{antoulas05}, the system $\Sigma = \{A,B,C\}$ can be embedded into a symmetric system $\hat{\Sigma} = \{\hat{A},\hat{B},\hat{C}\}$.
Since for each square matrix $A \in \R^{N \times N}$ there exists a symmetrizer $J = J^T$ such that $AJ = JA^T$ and thus the embedding system is given by:
\begin{align*}
 \hat{A} &= A, \\
 \hat{B} &= \begin{pmatrix} JC^T & B\end{pmatrix}, \\
 \hat{C} &= \begin{pmatrix} C \\ B^T J^{-1} \end{pmatrix}.
\end{align*}
If the system matrix $A$ is symmetric, and thus $J~=~\1$, the embedding of the system simplifies to:
\begin{align*}
 \hat{A} &= A, \\
 \hat{B} &= \begin{pmatrix} C^T & B\end{pmatrix}, \\
 \hat{C} &= \begin{pmatrix} C \\ B^T\end{pmatrix}.
\end{align*}
Even though the number of inputs and outputs is increased the number of states remains the same as in the original system.

\section{Parameter Identification and Combined Reduction} \label{pr}
The concept of controllability and especially observability extends to parametrized systems by treating the parameters as additional states.
These parameter states are constant over time and are assigned the parameters value as initial states:
\begin{align*}
 \dot{\breve{x}} &= \begin{pmatrix} \dot{x} \\ \dot{\theta} \end{pmatrix} = \begin{pmatrix} f(x(t),u(t),\theta) \\ 0 \end{pmatrix}, \\
 y &= g(x(t),u(t),\theta), \\
 \breve{x}_0 &= \begin{pmatrix} x_0 \\ \theta \end{pmatrix}.
\end{align*}
This augmented system, used in \cite{himpe13a}, can now be subject to a similar method to the direct truncation of the cross gramian for state reduction.

The cross gramian of the augmented system yields the \textbf{joint gramian} introduced in \cite{himpe13a}:
\begin{align*}
 W_J = \begin{pmatrix} W_X & \vline & W_M \\ \hline 0 & \vline & 0 \end{pmatrix} \in \R^{(n+p) \times (n+p)},
\end{align*}
with its upper left block ($W_X$) being the usual cross gramian of the system. The identifiability information of the parameters is encoded in $W_M$.
The parameter related information is then extracted by the Schur-complement of the symmetric part of the joint gramian, resulting in the \textbf{cross-identifiability gramian} $W_{\ddot{I}}$:
\begin{align*}
 W_{\ddot{I}} := W_M^* (\frac{1}{2}(W_X + W_X^T))^{-1} W_M.
\end{align*}
A singular value decomposition of $W_{\ddot{I}}$, provides a projection of the parameters that are sorted by their importance:
\begin{align*}
 W_{\ddot{I}} &= UDV, \\
 \Rightarrow \tilde{\theta} &= V \theta.
\end{align*}
Based on this projection the parameters can be partitioned and reduced,
\begin{align*}
 \tilde{\theta} &= \begin{pmatrix} \tilde{\theta_1} \\ \tilde{\theta_2} \end{pmatrix}, \\
 \Rightarrow \|\tilde{\theta_1}\|_1 &\approx \|\tilde{\theta}\|_1.
\end{align*}

As described in \cite{himpe13a}, by employing a truncation of states based on the singular values of $W_X$ and parameters based on the singular values of $W_{\ddot{I}}$ enables the combined reduction.

\section{Empirical Gramians} \label{eg}
Empirical gramians were introduced in \cite{lall99} and are solely based on simulations of the underlying control system.
These simulations use perturbations in input $u$ and initial states $x_0$ which are averaged.  
The required perturbations are organized into sets allowing a systematic perturbation of input $E_u \times R_u \times Q_u$ and initial states $E_x \times R_x \times Q_x$:
\begin{align*}
 E_u &= \{ e_i \in \R^j ; \|e_i\| = 1 ; e_i e_{j \neq i} = 0; i=1,\ldots,m \}, \\
 E_x &= \{ f_i \in \R^n ; \|f_i\| = 1 ; f_i f_{j \neq i} = 0; i=1,\ldots,n \}, \\
 R_u &= \{ S_i \in \R^{j \times j} ; S_i^* S_i = \1 ; i = 1,\ldots,s \}, \\
 R_x &= \{ T_i \in \R^{n \times n} ; T_i^* T_i = \1 ; i = 1,\ldots,t \}, \\ 
 Q_u &= \{ c_i \in \R ; c_i > 0 ; i = 1,\ldots,q \}, \\ 
 Q_x &= \{ d_i \in \R ; d_i > 0 ; i = 1,\ldots,r \}.
\end{align*}
Now the empirical cross gramian can be defined as follows (taken from \cite{himpe13a}):

For sets $E_u$, $E_x$, $R_u$, $R_x$, $Q_u$, $Q_x$, input $\bar{u}$ during steady state $\bar{x}$ with output $\bar{y}$, 
the \textbf{empirical cross gramian} $\hat{W}_X$ relating the states $x^{hij}$ of input $u^{hij}(t) = c_h S_i e_j + \bar{u}$ to output $y^{kla}$ of $x_0^{kla} = d_k T_l f_a + \bar{x}$, is given by:
 \begin{align*}
  \hat{W}_X &= \frac{1}{|Q_u| |R_u| m |Q_x| |R_x|} \sum_{h=1}^{|Q_u|} \sum_{i=1}^{|R_u|} \sum_{j=1}^m \sum_{k=1}^{|Q_x|} \sum_{l=1}^{|R_x|} \\ &\quad \quad \quad \cdot \frac{1}{c_h d_k} T_l \int_0^\infty \Psi^{hijkl}(t) d t \; T_l^*, \\
  &\Psi_{ab}^{hijkl}(t)  = f_b^* T_k^* \Delta x^{hij}(t) e_i^* S_h^* \Delta y^{kla}(t), \\
  &\Delta x^{hij}(t) = (x^{hij}(t)-\bar{x} ), \\
  &\Delta y^{kla}(t) = (y^{kla}(t)-\bar{y} ).
 \end{align*}

The joint gramian \cite{himpe13a} encapsulates the cross gramian, hence the \textbf{empirical joint gramian} is computed in the same manner as the empirical cross gramian, yet of the augmented system.
As shown in \cite{condon04}, the empirical gramians extend to time-varying systems, and thus can be applied in this setting for the hyperbolic networks.

\section{Uncertainty Quantification} \label{uq}
The connections between the network nodes, which are modeled by the components of the system matrix $A$, might contain uncertainties.
Due to the computation of empirical gramians based on simulations, potential uncertainties in initial state and external input can be incorporated by enlarging the corresponding set of perturbations respectively.
Hence, for an augmented system uncertainties in the parameters can also be included.
This allows robust model reduction.
Additionally, the parameter reducing projection can also be used to reduce, for example in a Gaussian setting, mean and covariance of a parameter distribution.

\section{Numerical Results} \label{nr}
To demonstrate the capabilities of this approach a synthetic hyperbolic network is utilized.
As described in section~\ref{hn}, the time varying system is growing with each time step.
This network with a maximum of $64$ nodes, thus a state dimension of $x \in \R^{64}$, and $8$ inputs and outputs is selected.
Furthermore, it is assumed that each connection is reciprocal, hence $A = A^T$ and $\theta \subset \R^{2016}$.
All possible connections of all nodes are treated as (time-varying) parameters in this setting,
\begin{align*}
 \dot{x} &= A(\theta(t)) x + Bu, \\
       y &= C x.
\end{align*}
Yet input matrix $B \in \R^{64 \times 8}$ and output matrix $C \in \R^{8 \times 64}$ are random and notably $C \neq B^T$, which requires an embedding into a symmetric system (see Section \ref{sr}).

First, in an offline phase, that has to be performed only once, a reduced order model is created.
The reduction procedure uses the empirical joint gramian of Section \ref{eg} that computes the cross gramian of the embedded augmented system:
\begin{align*}
 \dot{x} &= \begin{pmatrix} A(\theta(t)) \\ 0 \end{pmatrix} x + \begin{pmatrix} C^T & B \\ \quad 0\end{pmatrix} u, \\
       y &= \begin{pmatrix} C \\ B^T\end{pmatrix} x.
\end{align*}
Then, a reduction of states, based on the singular values of $W_X$, and of parameters, based on the singular values of $W_{\ddot{I}}$, is performed.
Second, in the online phase, the reduced model can be evaluated.

The computations are performed using the empirical gramian framework\footnote{See \url{http://gramian.de}} - \textbf{emgr} described in \cite{himpe13}.
Source code for the following experiments can be found at \mbox{\url{http://j.mp/ecc14_code}}.

\begin{figure}[ht!]
 \includegraphics[scale=0.4]{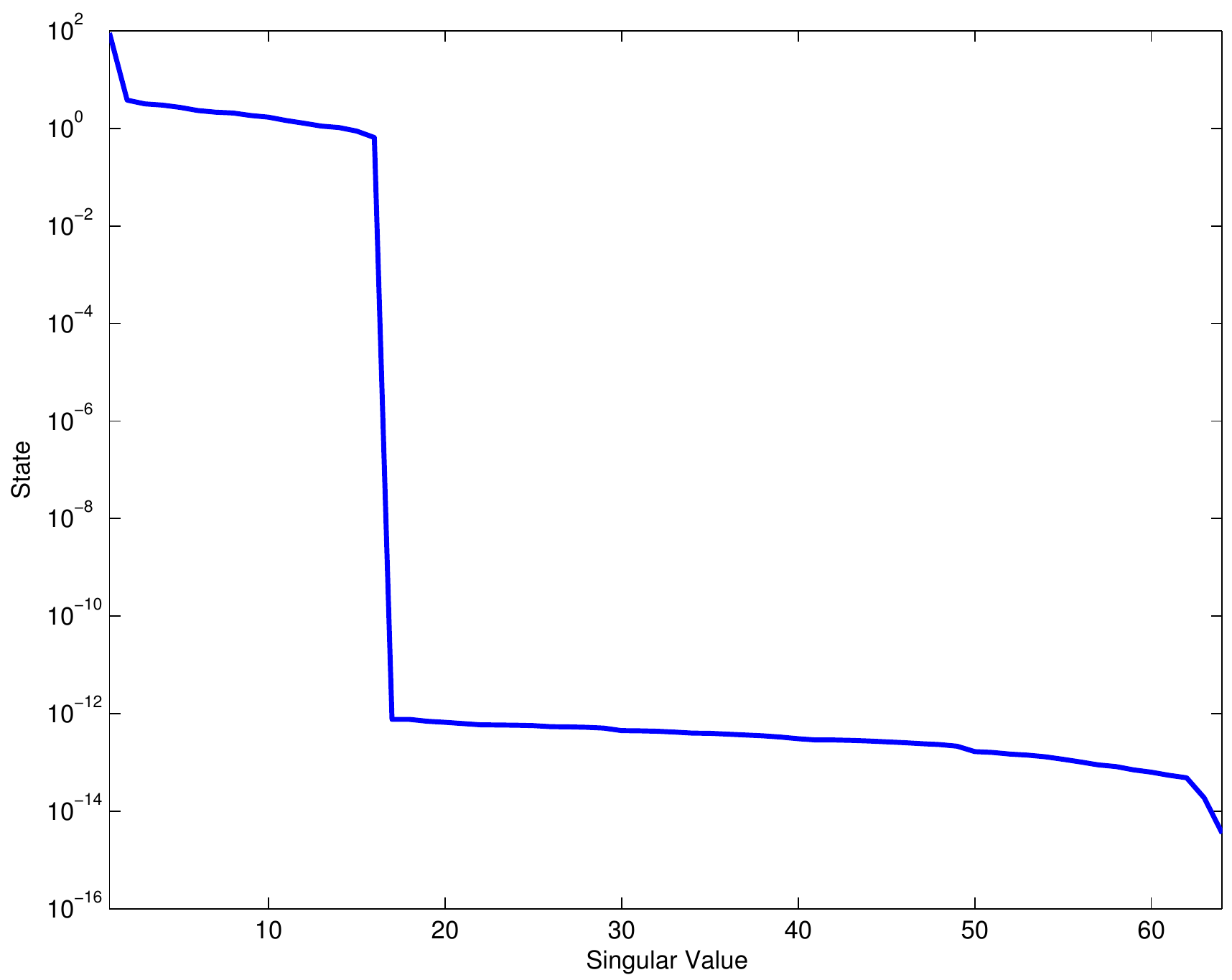}
 \caption{Distribution of the singular values of the cross gramian}
 \label{sv0}
 \includegraphics[scale=0.4]{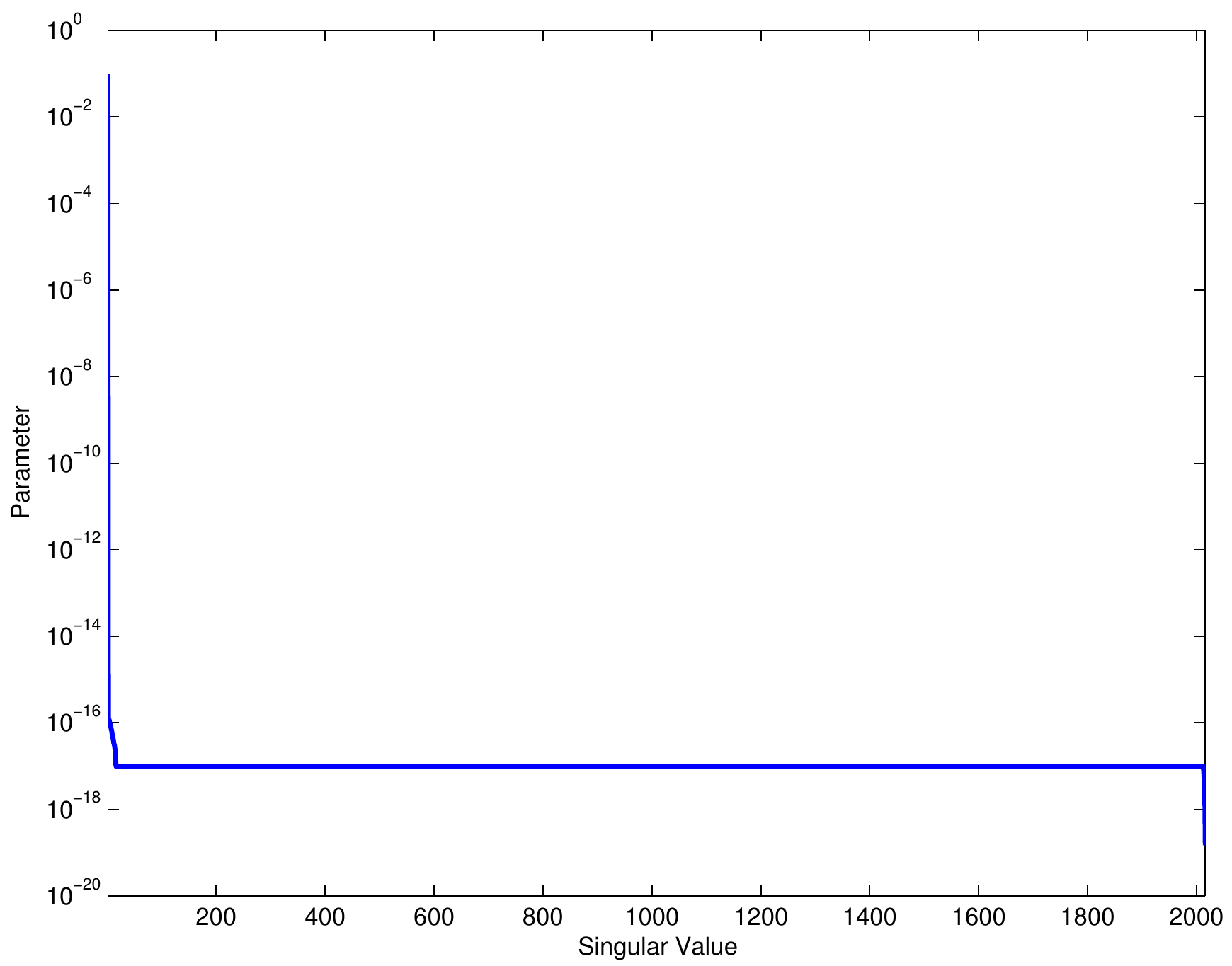}
 \caption{Distribution of the singular values of the cross identifiability gramian}
 \label{sv1}
\end{figure}

\begin{figure}[ht!]
 \includegraphics[scale=0.4]{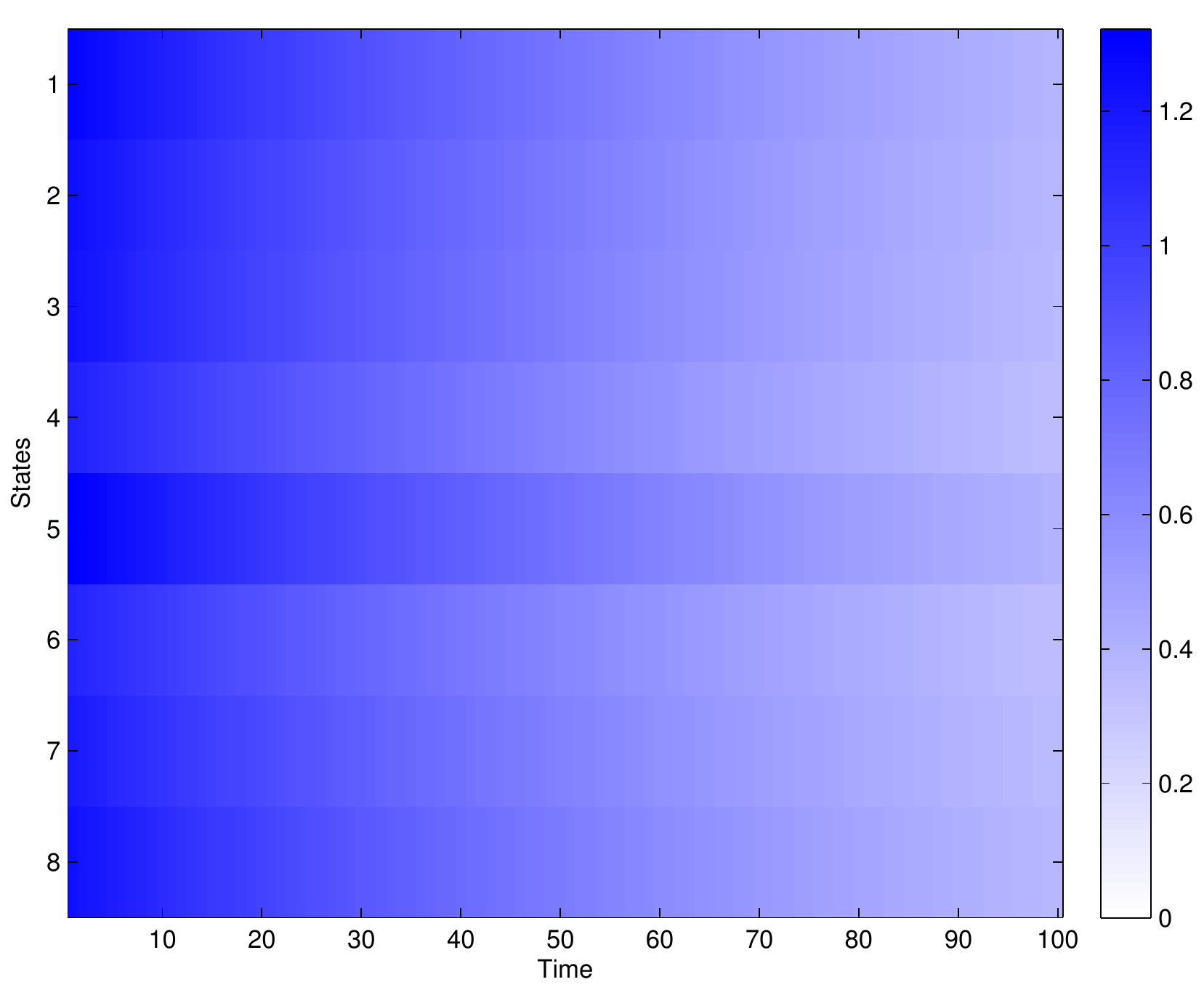}
 \caption{Impulse response of the full order model.}
 \label{hn0}
 \includegraphics[scale=0.4]{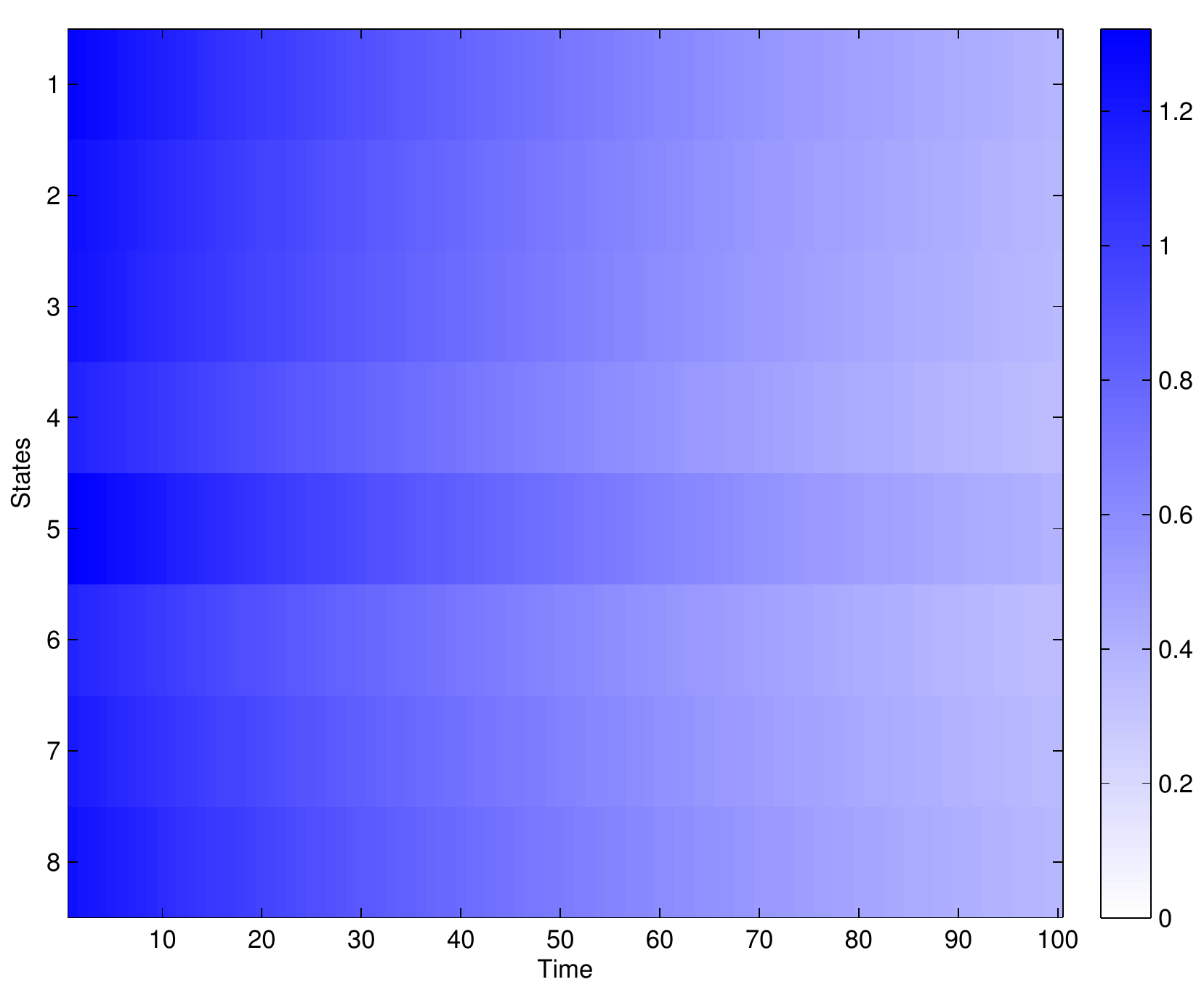}
 \caption{Impulse response of the reduced order model.}
 \label{hn1}
 \includegraphics[scale=0.4]{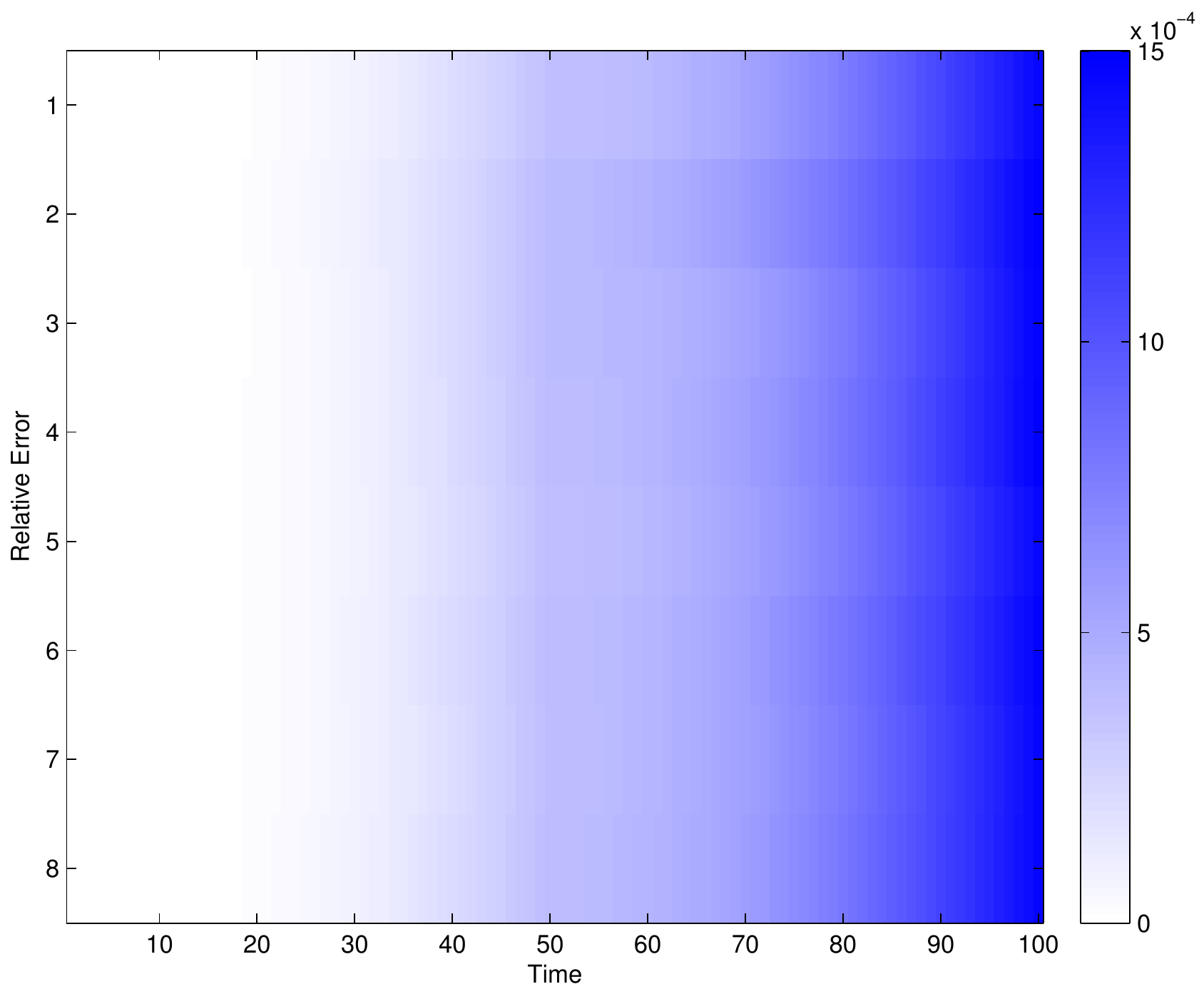}
 \caption{Relative error between the full and reduced order model.}
 \label{hn2}
\end{figure}

Since the network evolves during its evaluation the reduction has to be performed for unknown connectivity.
A distribution of the singular values of the empirical cross gramian $W_X$ and the empirical cross-identifiability gramian $W_{\ddot{I}}$ is given in Figure~\ref{sv0} and Figure~\ref{sv1}.
The singular values of these two empirical gramians describe the energy contained in the state and parameter respectively.
A reduction of the parameter space from dimension $2016$ to $65$ is suggested by the singular values of the cross-identifiability gramian.
For the state space a reduction from dimension $64$ to $19$ is performed, based on the singular values of the cross gramian.
The reduced model can be evaluated and compared to the full order model.
Figures~\ref{hn0} and \ref{hn1} show the impulse response of the full order and reduced order network, while the relative error between them is shown in Figure~\ref{hn2}.
Between the two time series of the full and reduced order model impulse response, the relative $L2$-error is $0.036\%$.
The sharp drop in the singular values of the cross-identifiability gramian determines the reduced order model dimension; a truncation of more parameter space dimensions will introduce a significant higher error.
The descent of singular values of the cross gramian also allows a graduated increase of error when truncating states.
To scan various location of the parameter space, only the low-dimensional reduced parameter space has to be scanned.

Comparing the reduced order model with the full order model, the combined reduction decreased the integration time by 21\% and memory requirements by 70\%.

\begin{table}[!h]
 \begin{tabular}{l|l|l|l}
  Original & Offline  & Online   & Relative \\
  Time (s) & Time (s) & Time (s) & L2-Error \\ \hline
  0.0203   & 277.1988 & 0.0160   & 0.00036  
 \end{tabular}
 \caption{Original Time, Offline Time, Online Time and Relative L2-Error; averaged over 100 simulations.}
\end{table}

\section{Conclusion}
The numerical experiments suggest that the empirical joint gramian \cite{himpe13a}, which is based on the empirical cross gramian, can be applied to reduce this type of control system with hyperbolic network structure.  
As shown, a linear time-varying control system in which also the parameter values vary over time, can be handled by the empirical gramians.
This concurrent reduction of state and parameter spaces enables, for example, scenarios in which the reduced model can be used to scan the parameter space. 

Using the (empirical) cross gramian of the embedded system is efficient, since the number of inputs and outputs is small compared to the number of states and no symmetrizer needs to be computed and inverted, $A$ has been chosen to be symmetric.
Yet, for networks with a non-symmetric system matrix $A$, a possibly costly computation of the symmetrizer and its inverse is required.
Thus, a generalization of the cross gramian to non-symmetric systems should be explored.



\begin{thebibliography}{99}

\bibitem{krioukov12}
D.~Krioukov, M.~Kitsak, R.~Sinkovits, D.~Rideout, D.~Meyer and M.~Bogu{\~n}{\'a}.6
\newblock Network Cosmology.
\newblock {\em Nature Publishing Group, Scientific reports}, 2:2012.

\bibitem{lall99}
S.~Lall, J.E.~Marsden, and S.~Glavaski.
\newblock Empirical model reduction of controlled nonlinear systems.
\newblock {\em Proceedings of the IFAC World Congress}, F:473--478, 1999.

\bibitem{condon04}
M.~Condon and R.~Ivanov.
\newblock Empirical balanced truncation of nonlinear systems.
\newblock {\em Journal of Nonlinear Science}, 14(5):405--414, 2004.

\bibitem{antoulas05}
A.C.~Antoulas.
\newblock Approximation of large-scale dynamical systems. 
\newblock {\em Society for Industrial Mathematics}, volume~6, 2005.

\bibitem{himpe13}
C.~Himpe and M.~Ohberger.
\newblock A Unified Famework for Empirical Gramians.
\newblock {\em Hindawi Journal of Mathematics}, 2013.

\bibitem{himpe13a}
C.~Himpe and M.~Ohlberger.
\newblock Cross-Gramian Based Combined State and Parameter Reduction for Large-Scale Control Systems.
\newblock {\em Submitted} Preprint at arXiv(math.OC): \href{http://arxiv.org/pdf/1302.0634.pdf}{1302.0634}, 2013.

\end{thebibliography}
\end{document}